\documentclass{article}
\usepackage{amssymb,latexsym,amscd}

\newcommand{\se}[1]{{\section{#1}} {\setcounter{equation}{0}}}
\newtheorem{theorem}{Theorem}[section]
\newtheorem{lm}{Lemma}[section]
\newtheorem{prop}{Proposition}[section]

\def\k{{K\"{a}hler }}
\def\ke{{K\"{a}hler-Einstein }}

\def\wp{{Weil-Peterson }}
\def\ma{{Monge-Amp\`{e}re }}
\begin{document}
\hbadness=10000
\title{{\bf Degeneration of \ke Manifolds II:}\\
{\Large {\bf The Toroidal Case}}}
\author{Wei-Dong Ruan\\
Department of Mathematics\\
University of Illinois at Chicago\\
Chicago, IL 60607\\}
\date{Revised: March 2004}
\footnotetext{Partially supported by NSF Grant DMS-0104150.}
\maketitle
%\tableofcontents
\begin{abstract}
In this paper we prove that the \ke metrics for a toroidal canonical degeneration family of \k manifolds with ample canonical bundles Gromov-Hausdorff converge to the complete \ke metric on the smooth part of the central fiber when the base locus of the degeneration family is empty. We also prove the incompleteness of the Weil-Peterson metric in this case. 
\end{abstract}
\se{Introduction}
This paper is a sequel of \cite{ruan3}. In algebraic geometry, when discussing the compactification of the moduli space of complex manifold $X$ with ample canonical bundle $K_X$, it is necessary to consider holomorphic degeneration family $\pi: {\cal X} \rightarrow B$, where $X_t = \pi^{-1}(t)$ are smooth for $t\not=0$, ${\cal X}$ and $X_0$ are $\mathbb{Q}$-Gorenstein, such that the canonical bundle of $X_t$ for $t\not=0$ and the dualizing sheaf of $X_0$ are ample. We will call such degeneration {\bf canonical degeneration}. We are interested in studying the degeneration behavior of the family of \ke metrics $g_t$ on $X_t$ when $t$ approaches 0. Following his seminal proof of Calabi conjecture (\cite{Yau2}), Yau (\cite{yau}) initiated the program of studying the application of \ke metrics to algebraic geometry with the belief that the behavior of \ke metrics should reflect the topological, geometric and algebraic structure of the underlying complex algebraic manifolds. According to this philosophy, one would expect the metric degeneration of the \ke manifolds to be closely related to the algebraic degeneration of the underlying algebraic manifolds. In \cite{Tian1}, Tian made the first important contribution along this direction. He proved that the \ke metrics on $X_t$ converge to the complete Cheng-Yau \ke metric on the smooth part of $X_0$ in the sense of Cheeger-Gromov, when ${\cal X}$ is smooth and the central fibre $X_0$ is the union of smooth normal crossing divisors $D_1,\cdots, D_L$, with a technical restriction that no three divisors have common intersection. Following the general framework in Tian's paper (\cite{Tian1}), \cite{Lu} and later \cite{ruan3} studied the general normal crossing case and removed the technical restriction in \cite{Tian1}. \\

In this paper, we generalize the result in \cite{ruan3} to the case when the central fibre $X_0$ is a union of toroidal orbifolds that results from the so-called toroidal canonical degeneration of smooth $X_t$ (see section 2 for definitions). The total space ${\cal X}$ for this kind of degeneration will be toroidal and generally not smooth. Please note that a toroidal canonical degeneration, where $X_0$ is not normal crossing in ${\cal X}$, can not be reduced to a normal crossing canonical degeneration. The normal crossing case is a very special case of toroidal canonical degeneration. For algebraic curve, a toroidal canonical degeneration is equivalent to Deligne-Mumford stable degeneration into stable curves.\\

In this paper, we always require an algebraic variety $X$ to possess a (set-theoretical) {\bf canonical} (Whitney) {\bf stratification} $\displaystyle X = \bigcup_{p\in \Sigma} D_p$ by smooth algebraic strata. By ``canonical" we mean that any other (Whitney) stratification $\displaystyle X = \bigcup_{p'\in \Sigma'} D'_{p'}$ by smooth strata is a refinement of the canonical (Whitney) stratification. More precisely, we have $D'_{p'} \subset D_p$ when $D'_{p'} \cap D_p \not= \emptyset$. For example, the toroidal varieties defined in section 2 satisfy such requirement.\\

The degeneration family is called {\bf base point free} if each smooth strata of $X_0$ is inside a smooth strata of ${\cal X}$. The smoothness condition of ${\cal X}$ when $X_0$ is normal crossing is equivalent to requiring the degeneration family to be base point free.  In some sense, toroidal canonical degenerations that we consider in this paper are generic base point free canonical degenerations. (Toroidal canonical degenerations and related concepts and constructions are discussed in section 2.)\\

Our first main theorem (proved in section 5) is the following.\\

\begin{theorem}
\label{aa}
Let $\pi: {\cal X} \rightarrow B$ be a toroidal canonical degeneration of \ke manifolds $\{X_t,g_{E,t}\}$ with ${\rm Ric}(g_{E,t}) = - g_{E,t}$. Then the \ke metrics $g_{E,t}$ on $X_t$ converge in the sense of Cheeger-Gromov to a complete Cheng-Yau \ke metric $g_{E,0}$ on the smooth part of the canonical limit $X'_0$ (which is a finite cover of the central fibre $X_0$).\\
\end{theorem}
To prove this theorem, we follow the three steps outlined in \cite{ruan3}. The first step is to construct certain smooth family of background \k metrics $\hat{g}_t$ on $X_t$ and their \k potential volume forms $\hat{V}_t$. The second step is to construct a smooth family of approximate \k metrics $g_t$ with \k form $\displaystyle\omega_t = \frac{i}{2\pi}\partial\bar{\partial}\log V_t$, where $V_t = h\hat{V}_t$ ($h$ is a function on ${\cal X}$) satisfies certain uniform estimate independent of $t$. The third step is to use Monge-Amp\`{e}re estimate of Aubin (\cite{A}) and Yau (\cite{Yau2}) to derive uniform estimate (independent of $t$) for the smooth family of \ke metrics $g_{E,t}$, starting with the smooth family of approximate \k metric $g_t$,  which is enough to ensure the Gromov-Hausdorff convergence of the family to the unique complete \ke metric $g_{E,0} = \{g_{0,i}\}_{i=1}^l$ on the smooth part of $X_0$. The first and the third steps are carried out in the very brief sections 3 and 5 and are virtually the same as in the normal crossing case \cite{ruan3}. The second step, carried out in section 4, is more involved than the simple global construction in \cite{ruan3}.\\ 

The following similar but much more non-trivial (comparing to \cite{ruan3}) estimate of the \wp metric near the degeneration, which implies the incompleteness of the \wp metric, is worked out in section 6.\\ 

\begin{theorem}
\label{ab}
The restriction of the \wp metric on the moduli space of complex structures to the toroidal canonical degeneration $\pi: {\cal X} \rightarrow B$ is bounded from above by a constant multiple of $\displaystyle\frac{dt\wedge d\bar{t}}{|\log |t||^3|t|^2}$. In particular, \wp metric is incomplete at $t=0$.\\
\end{theorem}
{\bf Note on notation:} We say $A\sim B$ if there exist constants $C_1,C_2>0$ such that $C_1B \leq A \leq C_2B$.\\

\se{Toroidal canonical degeneration}
In this section we introduce the concepts of toric degeneration and toroidal degeneration, and discuss the details of relevant stratification structures and the construction of compatible partition functions that we need for the construction of approximate metrics in section 4.\\
 
\subsection{Toric degenerations}
(Unless specified otherwise, the notations in this subsection will not be carried over to other parts of this paper.)\\ 

Let us first introduce the basic notions in toric geometry. An $(n+1)$-dimensional affine toric variety $A_{\sigma_0}$ is determined by a strongly convex $(n+1)$-dimensional integral polyhedral cone $\sigma_0$ in a rank $n+1$ lattice $\tilde{M}$. Let $\sigma_0(k)$ denote the set of $k$-dimensional subfaces of $\sigma_0$. Then $\sigma_0(n)$ corresponds to toric Weil divisors $\{D_i\}_{i\in \sigma_0(n)}$ in $A_{\sigma_0}$, and $\sigma_0(1)$ corresponds to toric Cartier divisors $\{(f_i)\}_{i\in \sigma_0(1)}$ in $A_{\sigma_0}$. For $i\in \sigma_0(1)$,

\[
(f_i) = \sum_{j\in \sigma_0(n)} a_{ij}D_j,
\]

where $a_{ij}$ is the natural pairing of the primitive elements in $i\in \sigma_0(1)$ and $j\in \sigma_0^\vee(1) \cong \sigma_0(n)$. $\sigma_0 ^{\vee}$ denotes the dual cone of $\sigma_0$.\\ 

A toric map $\pi: {\cal X} \rightarrow B \cong {\mathbb C}$ is called a {\bf toric degeneration}, if ${\cal X}=A_{\sigma_0}$ is an affine toric variety such that ${\cal X}\setminus X_0$ is the big open torus. Consequently $X_t$ for $t\not= 0$ are codimension one subtori in ${\cal X}\setminus X_0$. A toric degeneration is determined by a strongly convex integral polyhedral cone $\sigma_0\subset \tilde{M}$ with a marked primitive element $t$ in the interior of the cone $\sigma_0$. Under such notation, the central fibre is

\[
X_0 = D = \bigcup_{i\in \sigma_0(n)}D_i.
\]

Let $M = \tilde{M}/\mathbb{Z}\{t\}$. Since $t$ is in the interior of $\sigma_0$, the projection of $\sigma_0$ to $M$ determine a complete fan $\Sigma$ on $M$. Splittings $\tilde{M} \cong M\times \mathbb{Z}\{t\}$ can be parametrized (non-canonically) by $\mathbb{Z}$-valued linear functions on $M$. Each such splitting realizes $\sigma_0$ as a $\mathbb{Q}$-valued function $w_{\sigma_0}$ on the lattice $M$. In such a way, $\sigma_0$ can be understood as an equivalence class $[w_{\sigma_0}]$ (modulo $\mathbb{Z}$-valued linear functions) of convex piecewise linear $\mathbb{Q}$-valued functions on the lattice $M$ that are compatible with a complete fan $\Sigma$ in $M$. Let $\Sigma(k)$ denote the set of $k$-dimensional cones in $\Sigma$. Naturally $\Sigma(k) \cong \sigma_0(k)$ for $1\leq k \leq n$. We will use $\tilde{\sigma}\in \sigma_0(k)$ to denote the cone corresponding to $\sigma\in \Sigma(k)$.\\ 

For each $\sigma\in \Sigma$, there is an affine variety $A_\sigma = {\rm Spec}(\mathbb{C}[\sigma])$. For $\sigma,\sigma'\in \Sigma$ satisfying $\sigma\subset \sigma'$, there is a natural semi-group morphism $\sigma' \rightarrow \sigma$ that restricts to identity map on $\sigma \subset \sigma'$ and restricts to zero map on $\sigma'\setminus \sigma$, which induces the map $h_{\sigma\sigma'}: A_\sigma \rightarrow A_{\sigma'}$. Using $\{h_{\sigma\sigma'}\}_{\sigma,\sigma'\in \Sigma}$, we may glue the affine pieces $\{A_{\sigma}\}_{\sigma\in \Sigma}$ into the singular variety $X_\Sigma$. We have the following natural canonical (Whitney) stratification

\begin{equation}
\label{af}
X_\Sigma = \bigcup_{\sigma\in \Sigma} T_\sigma,\ {\rm where}\ T_\sigma = ({\rm Span}_{\mathbb{Z}}\sigma)^\vee\otimes_\mathbb{Z}\mathbb{C}^* = (M^\vee/\sigma^\perp)\otimes_\mathbb{Z}\mathbb{C}^*.
\end{equation}

In such a way, $\Sigma$ determines a singular variety $X_\Sigma$ that is mirror dual to the usual toric variety $P_\Sigma$ in certain sense.\\

For $\sigma\in \Sigma$, the natural injection $\tilde{\sigma} \hookrightarrow \sigma$ over $\mathbb{Z}$ induce a cover map $p_\sigma: A_\sigma \rightarrow A_{\tilde{\sigma}}$ and subsequently $q_\sigma = h_{\tilde{\sigma}\sigma_0}\circ p_\sigma: A_\sigma \rightarrow A_{\sigma_0}={\cal X}$. It is easy to check that $p_\sigma$, $q_\sigma$ for $\sigma\in \Sigma$ glue together to form the maps $p_\Sigma: X_\Sigma \rightarrow X_0$, $q_\Sigma: X_\Sigma \rightarrow {\cal X}$.\\  

Recall that a complex torus has a canonical {\bf toric holomorphic volume form}, and consequently a canonical real {\bf toric volume form}. Via this toric holomorphic volume form on the complex torus ${\cal X}\setminus X_0$, the dualizing sheaf $K_{\cal X}$ can be naturally identified with ${\cal O}_{\cal X}(-D)$. We call $\pi$ {\bf simple} when each divisor $D_i$ is of multiplicity one under $\pi$. Then the Cartier divisor $(t) = D$, and the dualizing sheaf $K_{\cal X}$ is a line bundle.\\
\begin{prop}
\label{ac}
A toric degeneration $\pi: {\cal X} \rightarrow B \cong {\mathbb C}$ is simple if and only if $[w_{\sigma_0}]$ is $\mathbb{Z}$-valued on $M$ if and only if $q_\Sigma: X_\Sigma \rightarrow {\cal X}$ is an imbedding (or equivalently, $p_\Sigma: X_\Sigma \rightarrow X_0$ is an isomorphism).
\end{prop}
{\bf Proof:} For $\tilde{\sigma}\in \sigma_0(n)$, it is straightforward to check that the multiplicity of $t$ along $D_{\tilde{\sigma}}$ is $|({\rm Span}_{\mathbb{Z}}\sigma)/({\rm Span}_{\mathbb{Z}}\tilde{\sigma})|$. Consequently, $\pi$ is simple if and only if for each $\sigma\in \Sigma(n)$, the natural injection $\tilde{\sigma} \hookrightarrow \sigma$ over $\mathbb{Z}$ is bijection (which amounts to that $[w_{\sigma_0}]$ is $\mathbb{Z}$-valued on $\sigma$) if and only if $p_\sigma: A_\sigma \rightarrow A_{\tilde{\sigma}}$ is an isomorphism for each $\sigma\in \Sigma(n)$. These local results together imply the proposition.
\hfill\rule{2.1mm}{2.1mm}\\
\begin{prop}
\label{ae}
For a toric degeneration $\pi: {\cal X} \rightarrow B \cong {\mathbb C}$, let $d$ be the smallest positive integer so that $d[w_{\sigma_0}]$ is $\mathbb{Z}$-valued. Then the canonical $d'$-fold base extension $\pi': {\cal X}' \rightarrow B'$ is a simple toric degeneration if and only if $d|d'$.
\end{prop}
{\bf Proof:} It is easy to see that the canonical $d'$-fold base extension $\pi': {\cal X}' \rightarrow B'$ is determined by $d'\sigma_0\subset \tilde{M}$. Since $[w_{d'\sigma_0}] = d'[w_{\sigma_0}]$ is $\mathbb{Z}$-valued if and only if $d|d'$, by proposition \ref{ac}, we get the desired conclusion.
\hfill\rule{2.1mm}{2.1mm}\\

{\bf Remark:} Propositions \ref{ac} and \ref{ae} are wellknown. (A special case of proposition \ref{ae}, where $\Sigma$ is simplicial fan, was proved and used in the proof of the semistable reduction theorem in \cite{Mu} by Mumford and Knudsen.) We provide simple proofs of them here for the convenience of the readers.\\

$\Sigma(1)$ can be equivalently interpreted as the set of primitive generating elements of 1-dimensional cones in $\Sigma$. The piecewise linear function $w_{\sigma_0}$ is determined by $\{w_m\}_{m\in \Sigma(1)}$, with $w_m = w_{\sigma_0}(m) \in \mathbb{Q}$ for $m\in \Sigma(1)$. The toric degeneration family can be equivalently characterized by the following family of toric immersions:

\[
i_t: N_{{\mathbb C}^*} \rightarrow {\mathbb C}^{|\Sigma(1)|}
\]

defined as $\{t^{w_m}z^m\}_{m\in \Sigma(1)}$, where $N = M^\vee$ and $N_{{\mathbb C}^*} = (N\otimes_{\mathbb{Z}}{\mathbb C})/N$. We are also interested in generalized toric degenerations, where $w_m \in \mathbb{R}$ are not necessarily rational.\\

{\bf Example:} The simplest toric degenerations that are not normal crossing are:

(1) $X_t = \{z\in \mathbb{C}^4|z_1z_2=z_3z_4=t\}$ (product of normal crossing degenerations).\\
(2) $X_t = \{z\in \mathbb{C}^4|z_1z_2=t,\ z_3z_4=tz_1\}$.
\hfill\rule{2.1mm}{2.1mm}\\

{\bf Remark}: {\it A priori}, the piecewise linear function $f$ generated by $\{w_m\}_{m\in \Sigma(1)}$ need not be convex. Then we may take the largest convex function $\tilde{f}\leq f$. The piecewise linear convex function $\tilde{f}$ will be generated by $\{w_m\}_{m\in \tilde{\Sigma}(1)}$, where $\tilde{\Sigma}(1)$ is a subset of $\Sigma(1)$. There is a natural projection $P: {\mathbb C}^{|\Sigma(1)|} \rightarrow {\mathbb C}^{|\tilde{\Sigma}(1)|}$. It is easy to check that $P$ induces an equivalence between the toric degeneration families determined by toric embeddings $i_t$ and $\tilde{i}_t = P\circ i_t$. Therefore, we only need to consider the case when $f$ is convex. For $f$ {\bf generic}, the fan it determines is a {\bf simplicial} fan. Namely the toric divisors are all {\bf toric orbifolds}.\\

\subsection{Toroidal degenerations}
A holomorphic degeneration $\pi: {\cal X} \rightarrow B = \{t\in \mathbb{C}:|t|<1\}$ is called a {\bf toroidal degeneration} if it is locally toric. Let

\begin{equation}
\label{ag}
X_0 = \bigcup_{k=0}^n X_0^{(k)} = D =\bigcup_{p\in \Sigma}D_p,\ X_0^{(k)} = \bigcup_{p\in \Sigma(k)}D_p,\ \Sigma = \bigcup_{k=0}^n \Sigma(k)
\end{equation}

be the canonical stratification for $X_0$, with $\{D_p\}_{p\in \Sigma}$ parametrizing all the strata and $X_0^{(k)}$ denoting the union of all $k$-dimensional strata. $\pi$ is called {\bf simple} if each divisor $\bar{D}_p$ is of multiplicity 1 under $\pi$ for $p\in \Sigma(n)$. Propositions \ref{ac} and \ref{ae} imply the following generalization to toroidal case.\\
\begin{prop}
\label{ah}
For a toroidal degeneration $\pi: {\cal X} \rightarrow B$, there exists an integer $d>0$ such that the canonical $d'$-fold base extension $\pi': {\cal X}' \rightarrow B'$ is a simple toroidal degeneration if and only if $d|d'$. $X'_0$ (which will be called the {\bf canonical limit}) is independent of $d'$ satisfying $d|d'$ with the natural finite cover map $X'_0\rightarrow X_0$.
\end{prop}
{\bf Proof:} Since $d'$-fold base extension is canonical and local, the $d>0$ here can be taken to be the lowest common multiple of the $d$'s specified in proposition \ref{ae} for all local toric models. When $d|d'$, namely $\pi'$ is simple, proposition \ref{ac} implies that $X'_0$ restricted to each local toric model can be identified with $X_\Sigma$ in proposition \ref{ac}, therefore is canonical and independent of $d'$.
\hfill\rule{2.1mm}{2.1mm}\\

Through proposition \ref{ah}, all discussions for toroidal degeneration can be reduced to discussions for simple toroidal degeneration via base extension. For this reason, we will always assume that $\pi$ is simple. Consequently, the dualizing sheaf $K_{\cal X}$ is a line bundle. For generic such degeneration $\pi$, the Weil divisors $\bar{D}_p$ for $p\in \Sigma(n)$ are {\bf toroidal orbifolds}. Without loss of generality and for simplicity of notations, we will also assume that each $\bar{D}_p$ does not self-intersect. \\

Choose a suitable tubular neighborhood $\tilde{U}_p$ of $\bar{D}_p$ for each $p\in \Sigma$ such that for any $p_1,p_2\in \Sigma$, we have

\[
\tilde{U}_{p_1}\cap\tilde{U}_{p_2} \subset \bigcup_{q\in\Sigma,D_q \in \bar{D}_{p_1}\cap \bar{D}_{p_2}} \tilde{U}_q.
\]

For each $p\in \Sigma$, we can construct a tubular neighborhood $U_p$ of $D_p$ as $\tilde{U}_p$ minus the union of divisors $\bar{D}_q$ for $q\in \Sigma(n)$ satisfying $D_p\not\subset \bar{D}_q$. We will also need $U_p^0 \subset U_p$ defined as $\tilde{U}_p$ minus the union of (slightly shrunk) $\overline{\tilde{U}}_q$ for $q\in \Sigma(n)$ satisfying $D_p\not\subset \bar{D}_q$. Let $D_p^0 = D_p \cap U_p^0$. Since

\[
\bigcup_{p\in \Sigma} U_p = \bigcup_{p\in \Sigma} U_p^0
\]

forms a neighborhood of $X_0$ that contains $X_t$ for $t$ small, many of our discussions on $X_t$ can be reduced locally to either $U_p\cap X_t$ or $U_p^0\cap X_t$ for $p\in \Sigma$. Notice that for any $p_1,p_2\in \Sigma$, we also have

\[
U_{p_1}\cap U_{p_2} \subset \bigcup_{q\in\Sigma,D_q \in \bar{D}_{p_1}\cap \bar{D}_{p_2}} U_q,\ \ U_{p_1}^0\cap U_{p_2}^0 \subset \bigcup_{q\in\Sigma,D_q \in \bar{D}_{p_1}\cap \bar{D}_{p_2}} U_q^0.
\]

Locally, $U_p = A_p\times D_p$ and $U_p^0 = A_p\times D_p^0$. $A_p$ is a neighborhood of the origin of the affine toric local model determined by the fan $\Sigma_p$ and the integral convex function $\{w_m\}_{m\in \Sigma_p(1)}$ (notation as in 2.1). Let $|p| := \dim D_p$ and $l = n-|p|$. $\Sigma_p(l)$ (which can be naturally identified with a subset of $\Sigma(n)$) corresponds to toroidal Weil divisors $\{D_q\cap U_p\}_{q\in \Sigma_p(l)\subset\Sigma(n)}$ in $U_p$. $\Sigma_p(1)$ corresponds to toroidal Cartier divisors $\{(s_m)\}_{m\in \Sigma_p(1)}$ in $U_p$ containing $D_p$. We may choose local coordinate $(t,z,\tilde{z})$ for $U_p$, $z=(z_1,\cdots,z_l)$, $\tilde{z}=(z_{l+1},\cdots,z_n)$, so that $s_m = t^{w_m}z^m$, $(t,z)$ and $\tilde{z}$ form coordinates for $A_p$ and $D_p$. $(z,\tilde{z})$ can be considered as coordinate for $X_t \cap U_p$. For $m\in \Sigma_p(1)$, $s_m$ can be viewed as a section of a line bundle on $U_p$ that defines the Cartier divisor. One can choose a Hermitian metric $\|\cdot\|_m$ on the line bundle over $U_p^0$ such that $\|s_m\|_m \leq 1$ and $\|s_m\|_m = 1$ outside a small neighborhood of the Cartier divisor $(s_m)$. More precisely, we require that $\|s_m\|_m = 1$ on $U_q^0$ for $q \in \Sigma$ when $s_m$ is nonvanishing on $D_q$.\\

For $p,q \in \Sigma$ satisfying $D_q \subset \bar{D}_p$, Cartier divisors in $U_p$ can be naturally extended to certain $\mathbb{Q}$-Cartier divisors in $U_q$ that can be expressed by the natural injective map $e_{pq}: \Sigma_p(1) \rightarrow \Sigma_q(1)$. By suitably adjusting the Hermitian metric of the line bundle, for $m\in \Sigma_p(1)$, we may assume that $\|s_m\|_m = \|s_{e_{pq}(m)}\|_{e_{pq}(m)}$ in the common domain $U_p^0\cap U_q^0$. It is easy to check that $e_{pq'} = e_{qq'}\circ e_{pq}$ for $q'\in \Sigma$ satisfying $D_{q'} \subset \bar{D}_q$. Therefore the Cartier divisor $(s_m)$ in $U_p^0$ for $m\in \Sigma_p(1)$ naturally extends to the $\mathbb{Q}$-Cartier divisor (still denoted by $(s_m)$) in $\tilde{U}_p$. $\|s_m\|_m$ for $m\in \Sigma_p(1)$ can similarly be extended from $U_p^0$ to $\tilde{U}_p$.\\

Let $\Sigma^p(1)$ denote the set of $q\in \Sigma(|p|+1)$ satisfying $D_q \subset \bar{D}_p$. For $q\in\Sigma^p(1)$, $D_p$ can be naturally identified with an element $[D_p]\in\Sigma_q(|q|-|p|)=\Sigma_q(1)$, which can also be viewed as a Cartier divisor $s_q$ in $\tilde{U}_p$ supported in $\tilde{U}_p\setminus U_p$. $\Sigma_p^p = \Sigma_p(1)\cup\Sigma^p(1)$ (resp. $\Sigma^p(1)$) can be characterized as the set of Cartier divisors on $\tilde{U}_p$ whose defining functions are not identically zero (resp. nowhere zero) on $X_0 \cap U_p$.\\

A {\bf (holomorphic) volume form} on $U_p\setminus D$ is called {\bf toroidal} if its pullback to the local toric model differ from the standard toric (holomorphic) volume form by a bounded nowhere zero (holomorphic) factor on $U_p$. By examining the holomorphic toric volume form, it is easy to see that a holomorphic toroidal volume form on ${\cal X}\setminus D$ can be naturally identified with a nowhere zero holomorphic section of $K_{\cal X}(D)$, or in another word, a meromorphic section of $K_{\cal X}$ with a pole of order 1 along $D$.\\

\subsection{Partition functions}
Let $\mu(x)$ be a smooth increasing function on $\mathbb{R}$ with bounded derivatives satisfying $\mu(x)=0$ for $x\leq 0$ and $\mu(x)=1$ for $x\geq 1$. Let $\min'(x_1,\cdots,x_l)$ be a smooth function with bounded derivatives that coincides with $\min(x_1,\cdots,x_l)$ when $\displaystyle\min_{i\not=j}(|x_i-x_j|)\geq 1$. (In another word, $\min'(x_1,\cdots,x_l)$ is a smoothing of $\min(x_1,\cdots,x_l)$ with bounded derivatives.)\\ 

For each $p\in \Sigma$ and $\eta>0$ large, we may define the smooth function

\[
\tilde{\mu}_p = \mu\left(\frac{1}{\log(\tau/\eta^2)}\textstyle\min'\left(\left\{\log(a_m/\eta)\right\}_{m\in \Sigma_p(1)}, \left\{\log(\tau/a_m\eta)\right\}_{m\in \Sigma^p(1)}\right)\right),
\]

where $\tau = -\log|t|^2$ and $a_m = \eta -\log\|s_m\|_m^2$. These will give us the partition functions $\{\mu_p\}_{p\in \Sigma}$, where $\displaystyle\mu_p = \tilde{\mu}_p\left(\sum_{p\in \Sigma} \tilde{\mu}_p\right)^{-1}$. We generally have $D_p^0 \subset {\rm supp}(\mu_p) \subset U_p$. The condition on $\|\cdot\|_m$ implies that 

\begin{equation}
\label{ad}
U_p^0 \cap {\rm supp}(\mu_q)=\emptyset\ \ {\rm when}\ \ D_p\not\subset \bar{D}_q.
\end{equation}

\se{Construction of the background metric}
For construction in this section to work, it is necessary to assume that the dualizing line bundle $K_{\cal X}$ of the total space ${\cal X}$ exists and is ample, which is valid in our situation. (The construction in this section is partially inspired by our work (\cite{ruan2}) on Bergmann metrics.) Recall that $K_{{\cal X}/B} = K_{\cal X}\otimes K_B^{-1}$ and $K_{X_t}= K_{{\cal X}/B}|_{X_t} \cong K_{\cal X}|_{X_t}$. (The last equivalence is not canonical, depending on the trivialization $K_B \cong {\cal O}_B$. We will use $dt$ to fix the trivialization of $K_B$.) Since $K_{X_t}$ is ample for all $t$, certain multiple $K^m_{X_t}$ will be very ample for all $t$. Equivalently, $K_{\cal X}^m$ is very ample on ${\cal X}$. It is not hard to find sections $\{\Omega_k\}_{k=0}^{N_m}$ of $K_{\cal X}^m$ that determine an embedding $e: {\cal X} \rightarrow \mathbb{CP}^{N_m}$, such that $\{\Omega_{t,k}\}_{k=0}^{N_m}$ forms a basis of $H^0(K_{X_t})$ for all $t$, where $\Omega_{t,k} = (\Omega_k\otimes (dt)^{-m})|_{X_t}$. $\{\Omega_{t,k}\}_{k=0}^{N_m}$ will determine a family of embedding $e_t: X_t \rightarrow \mathbb{CP}^{N_m}$ such that $e_t = e|_{X_t}$. Choose the Fubini-Study metric $\omega_{FS}$ on $\mathbb{CP}^{N_m}$, and define\\
\[
\hat{\omega} = \frac{1}{m}e^*\omega_{FS},\ \ \hat{\omega}_t = \hat{\omega}|_{X_t} = \frac{1}{m}e_t^*\omega_{FS}.
\]\\
Since $K_{\cal X}^m$ is very ample on ${\cal X}$, $\hat{\omega}$ is a smooth metric on ${\cal X}$. The \k potential of $\hat{\omega}$ and $\hat{\omega}_t$ are the logarithm of the volume forms\\
\[
\hat{V} = \left(\sum_{k=0}^{N_m} \Omega_{k}\otimes\bar{\Omega}_{k}\right)^{\frac{1}{m}},\ {\rm and}\ \hat{V}_t = \left(\sum_{k=0}^{N_m} \Omega_{t,k}\otimes\bar{\Omega}_{t,k}\right)^{\frac{1}{m}} = \left.\hat{V} \otimes (dt\otimes d\bar{t})^{-1}\right|_{X_t}.
\]

Since $K_{\cal X}^m$ is ample and therefore base point free, $\hat{V}$ is a non-degenerate smooth volume form on ${\cal X}$. Recall $(t) = D$. Hence $\displaystyle \frac{\hat{V}}{|t|^2}$ is a toroidal volume form on ${\cal X}$. On the other hand, $\displaystyle \frac{dt}{t}$ is the standard toric holomorphic form on $B$. Therefore,

\[
\hat{V}_t = \left(\sum_{k=0}^{N_m} \Omega_{t,k}\otimes\bar{\Omega}_{t,k}\right)^{\frac{1}{m}} = \left.\hat{V} \otimes (dt\otimes d\bar{t})^{-1}\right|_{X_t} = \left.\frac{\hat{V}}{|t|^2} \otimes \left(\frac{dt}{t}\otimes \frac{d\bar{t}}{\bar{t}}\right)^{-1}\right|_{X_t}
\]

is also toroidal, namely

\begin{equation}
\label{ba}
\hat{V}_t = \rho(t,z,\tilde{z})\left(\prod_{j=1}^l \frac{dz_jd\bar{z}_j}{|z_j|^2}\right)\left(\prod_{j=l+1}^n dz_jd\bar{z}_j\right)
\end{equation}

under the coordinate $(z,\tilde{z})$ for $X_t \cap U_p$, where $\rho(t,z,\tilde{z}) \sim 1$ is a smooth positive function on $U_p$.\\

Since $e: {\cal X} \rightarrow \mathbb{CP}^{N_m}$ is an embedding, locally in $U_p$, there exists a decomposition $e = \hat{e}\circ i_{\Sigma_p}$, where $i_{\Sigma_p} = (s_{\Sigma_p},\tilde{z}):{\cal X} \rightarrow \mathbb{C}^{|\Sigma_p(1)|+|p|}$ and $\hat{e}: \mathbb{C}^{|\Sigma_p(1)|+|p|} \rightarrow \mathbb{CP}^{N_m}$ are smooth embeddings and $s_{\Sigma_p} = (s_m)_{m\in \Sigma_p(1)}$. Therefore

\begin{equation}
\label{bb}
\hat{\omega} = \sum_{m,m'\in \Sigma_p(1)} g_{mm'}(s_{\Sigma_p},\tilde{z})ds_md\bar{s}_{m'} + ({\rm terms\ involving}\ d\tilde{z},d\bar{\tilde{z}}).
\end{equation}

\se{Construction of the approximate metric}
The approximate metric is constructed by gluing together appropriate metrics on the neighborhood of each strata by the partition functions constructed in 2.3.\\

For $p\in \Sigma$ and $m\in \Sigma_p^p$, in $U_p$, define

\[
h_p= \tau^{2(|\Sigma_p(1)|-l)}\prod_{m\in \Sigma_p^p}\frac{\eta^2}{a_m^2},\ \ a_m = \eta-\log\|s_m\|_m^2,\ \tau = -\log |t|^2.
\]

On $X_t$, let $V_t = h\hat{V}_t$, where $\displaystyle\log h = \sum_{p\in \Sigma}\mu_p\log h_p$, and let\\
\[
\omega_t = \frac{i}{2\pi}\partial\bar{\partial}\log V_t = \hat{\omega}_t + \frac{i}{2\pi}\partial\bar{\partial}\log h = \hat{\omega}_t + \gamma_t + \alpha_t,
\]

where

\[
\alpha_t = \sum_{p\in \Sigma}\mu_p\alpha_{t,p},\ \ \alpha_{t,p} = \frac{i}{\pi}\sum_{m\in \Sigma_p^p}\frac{1}{a_m^2}\partial a_m\bar{\partial}a_m,
\]
\[
\gamma_t = \sum_{p\in \Sigma}\mu_p\sum_{m\in \Sigma_p(1)}\frac{2}{a_m}{\rm Ric}(\|\cdot\|_m) + \frac{i}{2\pi}\sum_{p\in \Sigma}(\log h_p \partial\bar{\partial}\mu_p + \partial\log h_p \bar{\partial}\mu_p + \partial\mu_p \bar{\partial}\log h_p).
\]

The main result of this section is the estimate (Proposition \ref{cc} and Proposition \ref{cd}) on the approximate \k metric $g_t$ with the \k form $\omega_t$ on $X_t$.\\

Since $\Sigma_p$ is a simplicial fan, $\sigma\in \Sigma_p(l)$ naturally corresponds to a subset $S_\sigma \subset \Sigma_p(1)$ with $l$ elements.\\

\begin{prop}
\label{ce}
There exist $\lambda_1, \lambda_2 >0$ such that $\log\|s_m\|_m^2 \geq \lambda_2 \log |t|^2$ on $U_p^0$ for all $m\in \Sigma_p(1)$. And for any $x\in U_p^0$, $S_x = \{m\in \Sigma_p(1)|\log\|s_m(x)\|^2 \geq \lambda_1 \log |t|^2\} \subset S_\sigma$ for some $\sigma\in \Sigma_p(l)$.
\end{prop}
{\bf Proof:} Since $\log\|s_m\|_m^2 = \log|s_m|^2 + O(1)$, where $|s_m|$ is the absolute value of $s_m$ viewed as monomial in the toric local model, it is sufficient to prove the proposition for $\log|s_m|^2$ in place of $\log\|s_m\|_m^2$. For $m\in \Sigma_p(1)$, there exists $\sigma\in \Sigma_p(l)$ such that $-m$ belongs to the cone spanned by $S_\sigma$. Namely 

\[
m = -\sum_{m' \in S_\sigma} b_{m'} m',
\]

where $b_{m'}\geq 0$ for all $m' \in S_\sigma$. Therefore

\[
\log |s_m|^2 = w_m\log |t|^2 + \log|z^m|^2 = w_m\log |t|^2 - \sum_{m' \in S_\sigma} b_{m'}\log|z^{m'}|^2 
\]
\[
= (w_m + \sum_{m' \in S_\sigma} b_{m'}w_{m'})\log |t|^2 - \sum_{m' \in S_\sigma} b_{m'}\log |s_{m'}|^2 
\]
\[
\geq (w_m + \sum_{m' \in S_\sigma} b_{m'}w_{m'})\log |t|^2.
\]

We may take $\lambda_2$ to be the maximum of such $(w_m + \sum_{m' \in S_\sigma} b_{m'}w_{m'})$.\\

Take a subset $S'\subset S$ such that $S'$ span a simplicial cone and $S' \not\subset S_\sigma$ for any $\sigma\in \Sigma_p(l)$. There exists a linear function $v_m$ on $M$ such that $v_m = \frac{\log |s_m|^2}{\log |t|^2} \leq \lambda_1$ for $m\in S'$ and $|v_m| \leq C\lambda_1$ for $m\in \Sigma_p(1)\setminus S'$. Then

\[
w'_m = \frac{\log |s_m|^2}{\log |t|^2} - v_m = w_m + \frac{\log|z^m|^2}{\log |t|^2} - v_m
\]

is an adjustment of $w_m$ by a linear function on $M$, such that $w'_m =0$ for $m\in S'$. Since $S' \not\subset S_\sigma$ for any $\sigma\in \Sigma_p(l)$. The strict convexity of $\{w_m\}_{m\in \Sigma_p(1)}$ implies that there exists an $m'\in \Sigma_p(1)\setminus S'$ such that $w'_{m'}<0$ is the smallest. Take $\lambda_3$ to be the maximum of such $w'_{m'}<0$ for all possible $S'$, which have only finite many possibilities. Then $\lambda_3<0$ and

\[
\frac{\log |s_{m'}|^2}{\log |t|^2} = w'_{m'} + v_{m'} \leq \lambda_3 + C\lambda_1.
\]

We may take $\lambda_1>0$ to be small so that $\lambda_3 + C\lambda_1<0$. Then $|s_{m'}|^2$ has to be big, contradicting the fact that $|s_{m'}|^2$ is small in $U_p$. Therefore, $S\subset S_\sigma$ for some $\sigma\in \Sigma_p(l)$.
\hfill\rule{2.1mm}{2.1mm}\\
\begin{lm}
\label{ch}
\[
\gamma_t = O(\hat{\omega}_t/\eta) + O((\hat{\omega}_t + \alpha_t)/\log \tau),\ \ {\rm where}\ \ \tau = -\log |t|^2.
\]
\end{lm}
{\bf Proof:} In the argument of this paper, we will always first fix $\eta>0$ large and then take $\tau$ large according to the fixed $\eta$. By our construction, $a_m\geq \eta$ is large. Hence

\[
\sum_{p\in \Sigma}\mu_p\sum_{m\in \Sigma_p(1)}\frac{2}{a_m}{\rm Ric}(\|\cdot\|_m) = O(\hat{\omega}_t/\eta).
\]

For any $x\in X_t$, there exist a $q\in \Sigma$ such that $x\in X_t\cap U_q^0$. Since

\[
\sum_{p\in \Sigma} \mu_p =1,\ \ \sum_{p\in \Sigma} \partial\bar{\partial}\mu_p =0.
\]

We have

\[
\sum_{p\in \Sigma}\log h_p \partial\bar{\partial}\mu_p = \sum_{p\in \Sigma}(\log h_p - \log h_q) \partial\bar{\partial}\mu_p.
\]

Since $U_q^0 \cap {\rm supp}(\mu_p)=\emptyset$ when $D_q\not\subset \bar{D}_p$ according to (\ref{ad}), we may consider only those $p\in \Sigma$ satisfying $D_q\subset \bar{D}_p$. Then there are the natural inclusions $\Sigma^q(1) \subset\Sigma^p(1)$, $\Sigma_p(1) \subset\Sigma_q(1)$ and the Cartier divisors in $\Sigma^p(1) \setminus \Sigma^q(1)$ vanishing along $D_q$ can be naturally identified with a subset of $\Sigma_p(1) \setminus \Sigma_q(1)$. Under such identifications $\Sigma_q^q \cap\Sigma_p^p$ is defined. For any $m\in\Sigma_p^p \setminus \Sigma_q^q$, $(s_m)\cap \bar{D}_q=\emptyset$. Consequently, $\|s_m\|_m^2= 1$ and $a_m = \eta$ on $\tilde{U}_q$ for $m\in\Sigma_p^p \setminus \Sigma_q^q$. Hence

\[
\log h_p - \log h_q = 2\sum_{m\in \Sigma_q^q \setminus\Sigma_p^p}\log \frac{a_m}{\tau\eta}
\]

is bounded on ${\rm supp}(\mu_p)\cap U_q \subset U_p\cap U_q$. From the explicit expressions of $\mu_p$ and $h_p$, it is straightforward to check that $\partial\bar{\partial}\mu_p = O(1/\log \tau)$, $\partial\log h_p = O(1)$ and $\partial\mu_p = O(1/\log \tau)$ with respect to the Hermitian metric $\hat{\omega}_t + \alpha_t$. (Such kind of verification is more carefully done in the proof of proposition \ref{cb} using (\ref{ca}).) Consequently $\gamma_t = O(\hat{\omega}_t/\eta) + O((\hat{\omega}_t + \alpha_t)/\log \tau)$.
\hfill\rule{2.1mm}{2.1mm}\\

For $\sigma\in \Sigma_p(l)$, let $A_\sigma(x) = \displaystyle\min_{m\in \Sigma_p(1)\setminus S_\sigma} a_m(x)$ and

\[
U_{p\sigma}^0 = U_{p\sigma}\cap U_p^0,\ U_{p\sigma} = \{x\in U_p|A_\sigma(x) \geq A_{\sigma'}(x)\ {\rm for}\ \sigma'\in \Sigma_p(l)\}.
\]

Then the proposition \ref{ce} implies that $A_\sigma(x) \geq \lambda_1\tau >0$ for $x \in U_{p\sigma}^0$ and

\begin{prop}
\label{cg}
For $t$ small enough, we have

\[
X_t \cap U_p^0 = \bigcup_{\sigma\in \Sigma_p(l)} X_t \cap U_{p\sigma}^0,
\]

and $a_m^2 \sim (\log |t|^2)^2$ in $U_{p\sigma}^0$ for $m\in \Sigma(1)\setminus S_\sigma$.
\hfill\rule{2.1mm}{2.1mm}\\
\end{prop}
For $S_\sigma = \{m_1,\cdots,m_l\}$, on $U_{p\sigma}$, we may choose coordinate $z = \{z_k\}_{k=1}^l = \{s_{m_k}\}_{k=1}^l$. By adjusting the convex function $w=\{w_m\}_{m\in \Sigma_p(1)}$ by linear function, we may assume that $w_{m} = 0$ for $m\in S_\sigma$ and $w_m>0$ for $m\in \Sigma_p(1)\setminus S_\sigma$. Then we have $s_m = t^{w_m}z^m$, where $m = \{m^k\}_{k=1}^l$ also denotes the coordinate of $m$ with respect to the basis $\{m_k\}_{k=1}^l$. It is easy to see that this coordinate $z$ is a special case of the toroidal coordinate $z$ defined in section 2. Let

\[
\alpha_{p\sigma} = \frac{i}{\pi}\sum_{m\in S_\sigma}\frac{1}{a_m^2}\partial a_m\bar{\partial}a_m,\ \alpha_{t,p\sigma} = \alpha_{p\sigma}|_{X_t}.
\]

\begin{prop}
\label{cf}
\[
\alpha_{p\sigma} \leq \alpha \leq C(w)\alpha_{p\sigma}
\]
along $z$ direction in $U_{p\sigma}^0$. Consequently,
\[
C_1\left(\prod_{m\in S_\sigma}\frac{1}{a_m^2}\right)\hat{V}_t \leq \omega_t^n \leq C_2\left(\prod_{m\in S_\sigma}\frac{1}{a_m^2}\right)\hat{V}_t,\ \ {\rm in}\ \ U_{p\sigma}^0\cap X_t.
\]
\end{prop}
{\bf Proof:} By the definition of $U_{p\sigma}^0$, clearly $\alpha_{p\sigma} \leq \alpha \leq C(w)\alpha_{p\sigma}$ along $z$ direction in $U_{p\sigma}^0$. Therefore $\omega_t \sim  \hat{\omega}_t + \alpha_t \sim \hat{\omega}_t + \alpha_{t,p\sigma}$ according to lemma \ref{ch}. Since $\alpha_{t,p\sigma}^{l+1}=0$. In $X_t \cap U_{p\sigma}^0$, we have

\[
\omega_t^n \sim  (\hat{\omega}_t + \alpha_{t,p\sigma})^n \sim \hat{\omega}_t^{n-l} \wedge \alpha_{t,p\sigma}^l.
\]

According to formula (\ref{ba}),

\[
\hat{V}_t = \rho(z)\left(\prod_{j=1}^l \frac{dz_jd\bar{z}_j}{|z_j|^2}\right)\left(\prod_{j=l+1}^n dz_jd\bar{z}_j\right).
\]

Hence 

\[
\hat{\omega}_t^{n-l} \wedge \alpha_{t,p\sigma}^l \sim \left(\prod_{m\in S_j}\frac{\partial a_m\bar{\partial}a_m}{a_m^2}\right)\left(\prod_{j=l+1}^n dz_jd\bar{z}_j\right) \sim \left(\prod_{m\in S_p}\frac{1}{a_m^2}\right)\hat{V}_t.
\]
\hfill\rule{2.1mm}{2.1mm}\\

Notice that $V_t = h\hat{V}_t$ is the \k potential of $\omega_t$. Assume

\[
e^{-\phi_t} = \frac{\omega_t^n}{V_t}.
\]

\begin{prop}
\label{cb}
$|\phi_t|$ is bounded independent of $t$.
\end{prop}
{\bf Proof:} According to proposition \ref{cg}, it is sufficient to verify in each $U_{p\sigma}^0 \cap X_t$ for $p\in \Sigma$ and $\sigma\in \Sigma_p(l)$. Proposition \ref{cf} implies that

\[
\frac{\omega_t^n}{V_t} \sim \eta^{2|\Sigma_p(1)|} \sim 1\ \ {\rm in}\ U_{p\sigma}^0 \cap X_t.
\]

Therefore $|\phi_t|$ is bounded independent of $t$. 
\hfill\rule{2.1mm}{2.1mm}\\

Let $g_t$ denote the \k metric corresponding to the \k form $\omega_t$, then we have\\
\begin{prop}
\label{cc}
The curvature of $g_t$ and its derivatives are all uniformly bounded with respect to $t$.
\end{prop}
{\bf Proof:} On a Riemannian manifold $(M,g)$, we call a basis $\{v_i\}$ {\bf proper} if the corresponding metric matrix satisfies $C_1(\delta_{ij}) \leq (g_{ij}) \leq C_2(\delta_{ij})$ for $C_1,C_2>0$. To verify that the curvature of the Riemannian metric $g$ and all its covariant derivatives are bounded, it is sufficient to find a proper basis $\{v_i\}$ satisfying that the coefficients of $[v_i,v_j]$ with respect to the basis $\{v_i\}$ and all their derivatives with respect to $\{v_i\}$ are bounded, such that $g_{ij}$ and all their derivatives with respect to $\{v_i\}$ are bounded.\\

According to proposition \ref{cg}, it is sufficient to verify in each $U_{p\sigma}^0 \cap X_t$. Let $W_j = a_{m_j}z_j\displaystyle\frac{\partial}{\partial z_j}$ for $1\leq j \leq l$ and $W_j = \displaystyle\frac{\partial}{\partial z_j}$ for $l+1\leq j \leq n$. According to proposition \ref{cf}, it is straightforward to check that the basis $\{W_j,\bar{W}_j\}_{j=1}^n$ is proper in $U_{p\sigma}^0 \cap X_t$. Namely $C_1(\delta_{ij}) \leq (g_{i\bar{j}}) \leq C_2(\delta_{ij})$ for some $C_1,C_2>0$, where $(g_{i\bar{j}})$ denotes the metric matrix with respect to the basis $\{W_j,\bar{W}_j\}_{j=1}^n$. (For the upper bound estimate, we need $\displaystyle\frac{a_{m_j}}{a_m}$ to be bounded for $1\leq j \leq l$ and $m\in \Sigma_I(1)\setminus S_i$, which is due to our restriction to $U_{p\sigma}^0$.)\\ 

For $1\leq j \leq l$, $\|s_{m_j}\|_{m_j}^2 = \rho_j|z_j|^2.$

\[
W_k(a_{m_j}) = \frac{W_k(\|s_{m_j}\|_{m_j}^2)}{\|s_{m_j}\|_{m_j}^2} = \frac{W_k(\rho_j)}{\rho_j} + \frac{W_i(|z_j|^2)}{|z_j|^2}.
\]

\[
W_k(a_{m_j}) = a_{m_k}(z_k\frac{\partial\log \rho_j}{\partial z_k} + \delta_{kj})\ \ {\rm for}\ 1\leq k \leq l.
\]

\[
W_k(a_{m_j}) = \frac{\partial\log \rho_j}{\partial z_k}\ \ {\rm for}\ l+1\leq k \leq n.
\]

The functions

\begin{equation}
\label{ca}
\begin{array}{ll}
\displaystyle \frac{1}{a_{m_j}},\ s_mP(a),\ \bar{s}_mP(a),\ z_jP(a_{m_j}),\ \bar{z}_jP(a_{m_j}),\\
\displaystyle \frac{a_{m_j}}{a_m},\ \frac{\log|t|^2}{a_m},\ \frac{a_m}{\log|t|^2},\ \ {\rm for}\ m\in \Sigma_p(1)\setminus S_\sigma,\ 1\leq j \leq l.
\end{array}
\end{equation}

are all bounded in $U_{p\sigma}^0 \cap X_t$, where $P(a)$ is a polynomial on $(\{a_{m_j}\}_{j=1}^l,\log t)$ and $P(a_{m_j})$ is a polynomial on $a_{m_j}$. Above computations imply that the derivatives of functions in (\ref{ca}) with respect to $\{W_j,\bar{W}_j\}_{j=1}^n$ are smooth functions of terms in (\ref{ca}) and other smooth bounded terms. Therefore they are bounded.\\

It is straightforward to check that $g_{i\bar{j}}$ and the coefficients of $[W_j,W_k]$, $[W_j,\bar{W}_k]$, $[\bar{W}_j,\bar{W}_k]$ with respect to the basis $\{W_j,\bar{W}_j\}_{j=1}^n$ are all smooth functions of terms in (\ref{ca}) and other bounded smooth terms. Consequently, any derivatives of theirs with respect to $\{W_j,\bar{W}_j\}_{j=1}^n$ are also smooth functions of terms in (\ref{ca}) and other bounded terms, therefore, are all bounded.
\hfill\rule{2.1mm}{2.1mm}\\
\begin{prop}
\label{cd}
For any $k$, $\|\phi_t\|_{C^k,g_t}$ is uniformly bounded with respect to $t$.\\
\end{prop}
{\bf Proof:} Similar as in the proof of the previous proposition, in $U_{p\sigma}^0 \cap X_t$, it is straightforward to check according to proposition \ref{cb} and the explicit expression of $\phi_t$ that $\phi_t$ is a bounded smooth function of terms in (\ref{ca}) and other smooth bounded terms. Consequently, all multi-derivatives of $\phi_t$ with respect to $\{W_j,\bar{W}_j\}_{j=1}^n$ are smooth functions of terms in (\ref{ca}) and other smooth bounded terms. Therefore they are bounded.
\hfill\rule{2.1mm}{2.1mm}\\

\se{Construction of \ke metric via complex Monge-Amp\`{e}re}
In this section, we will use the same notions as in the previous sections. In \cite{Tian1}, using the Monge-Amp\`{e}re estimate of Aubin and Yau, Tian essentially proved the following.\\

\begin{theorem}
\label{db}
(Tian) Assume that $\phi_t$, the curvature of $g_t$ and their multi-derivatives are all bounded uniformly independent of $t$, then the \ke metric $g_{E,t}$ on $X_t$ will converge to the complete Cheng-Yau \ke metric $g_{E,0}$ on $X_0\setminus{\rm Sing}(X_0)$ in the sense of Cheeger-Gromov: there are an exhaustion of compact subsets $F_\beta \subset X_0\setminus{\rm Sing}(X_0)$ and diffeomorphisms $\psi_{\beta,t}$ from $F_\beta$ into $X_t$ satisfying:\\
(1) $\displaystyle X_t\setminus \bigcup_{\beta=1}^\infty \psi_{\beta,t}(F_\beta)$ consists of finite union of submanifolds of real codimension 1;\\
(2) for each fixed $\beta$, $\psi_{\beta,t}^*g_{E,t}$ converge to $g_{E,0}$ on $F_\beta$ in $C^k$-topology on the space of Riemannian metrics as $t$ goes to $0$ for any $k$.
\hfill\rule{2.1mm}{2.1mm}\\
\end{theorem}
{\bf Proof of theorem \ref{aa}:} Proposition \ref{ah} reduces the theorem to the case that $\pi: {\cal X} \rightarrow B$ is simple, which is a direct corollary of theorem \ref{db} and propositions \ref{cb}, \ref{cc}, \ref{cd}. 
\hfill\rule{2.1mm}{2.1mm}\\

It is easy to see that our construction actually implies the following asymptotic description of the family of \ke metrics.\\

\begin{theorem}
\label{da}
\ke metric $g_{E,t}$ on $X_t$ is uniformly quasi-isometric to the explicit approximate metric $g_t$. More precisely, there exist constants $C_1,C_2>0$ independent of $t$ such that $C_1 g_t \leq g_{E,t}\leq C_2 g_t$.\\
\end{theorem}
{\bf Proof:} 
The uniform $C^0$-estimate of the complex Monge-Amp\`{e}re equations implies that $C_1 \omega_t^n \leq \omega_{E,t}^n\leq C_2 \omega_t^n$ for some $C_1,C_2>0$. The uniform $C^2$-estimate of the complex Monge-Amp\`{e}re equations implies that ${\rm Tr}_{g_t} g_{E,t}$ is uniformly bounded from above. Combining these two estimates, we get our conclusion. 
\hfill\rule{2.1mm}{2.1mm}\\

\se{\wp metric near degeneration}
In this section, we will start with the discussion of the toric case, which is of independent interest and the estimate is more precise. Then we will proceed to the global toroidal case.\\

\subsection{the toric case}
{\bf Note}: The notations in this subsection are the same as in subsection 2.1. Unless specified otherwise, the notations in this section will not be carried over to other parts of this paper.\\ 

{\bf Example:} Consider a toric degeneration $\pi: {\cal X} \rightarrow B \cong {\mathbb C}$ determined by a complete fan $\Sigma$ in $M$ and an integral piecewise linear convex function determined by $\{w_m\}_{m\in \Sigma(1)}$. For $i\in \Sigma(n)$ assume $w_{m} = 0$ for $m\in S_i$ and $w_m>0$ for $m\in \Sigma(1)\setminus S_i$. With $S_i = \{m_1,\cdots,m_n\}$ and toric coordinate $z_j = s_{m_j}$ for $1\leq j \leq n$, we have

\[
\omega = \frac{i}{\pi}\sum_{j=1}^n\frac{dz_j\wedge d\bar{z}_j}{|z_j|^2(\log|z_j|^2)^2} + \frac{i}{\pi}\sum_{m\in \Sigma(1)\setminus S_i}\frac{ds_m\wedge d\bar{s}_m}{|s_m|^2(\log|s_m|^2)^2}.
\]

Let

\[
W= \frac{\nabla \log t}{|\nabla \log t|^2},
\]

then $\bar{\partial} W|_{X_t}$ is a natural representative of Kodaira-Spencer deformation class in the Dolbeaut cohomology $H^1(T_{X_t})$. $W$ can also be determined by the conditions $\displaystyle\pi_*W = t\frac{d}{dt}$ and $i(W)\omega|_{X_t} = 0$ for all $t$. Let $a_j = \log|z_j|^2$ and $a_m =\log|s_m|^2$ for $m\in \Sigma(1)\setminus S_i$. We will use $\rho = 1+ O(a_j/a_m)$ to denote a bounded smooth function on $a_j/a_m$ for $1\leq j \leq n,m\in \Sigma(1)\setminus S_i$. (Here $O(a_j/a_m)$ is a shorthand for $O(a_j/a_m,1\leq j \leq n,m\in \Sigma(1)\setminus S_i)$.) It is straightforward to derive that\\
\[
\omega_t = \omega|_{X_t} = \frac{i}{\pi}\sum_{j,k=1}^ng_{j\bar{k}}\frac{\partial a_j}{a_j}\wedge\frac{\bar{\partial}a_k}{a_k},
\]
\[
g_{j\bar{k}} = \delta_{jk} + a_ja_kO\left(\frac{1}{a_m^2}\right),\ \ g^{j\bar{k}} = \delta_{jk} + a_ja_kO\left(\frac{1}{a_m^2}\right),
\]
\[
\omega^n_t = n!\left(\frac{i}{\pi}\right)^n\rho \prod_{j=1}^n\frac{dz_j\wedge d\bar{z_j}}{a_j^2|z_j|^2} = n!\left(\frac{1}{\pi}\right)^n\rho \prod_{j=1}^n\frac{da_j\wedge d\theta_j}{a_j^2}.
\]
\hfill\rule{2.1mm}{2.1mm}\\
\begin{lm}
\label{eb}
\[
W = t\frac{\partial}{\partial t} - \sum_{j=1}^n\sum_{m\in \Sigma(1)\setminus S_i} w_mm^j\frac{a_j^2}{a_m^2}\rho z_j\frac{\partial}{\partial z_j}.
\]
\end{lm}
{\bf Proof:} Since $\displaystyle\pi_*W = t\frac{d}{dt}$, we may assume that $\displaystyle W = t\frac{\partial}{\partial t} + \sum_{j=1}^nq_j z_j\frac{\partial}{\partial z_j}$. $i(W)\omega|_{X_t} = 0$ implies that

\[
\sum_{m\in \Sigma(1)\setminus S_i} w_m\frac{d\bar{s}_m}{\bar{s}_ma_m^2} + \sum_{j=1}^nq_j \frac{d\bar{z}_j}{\bar{z}_ja_j^2} + \sum_{m\in \Sigma(1)\setminus S_i}\sum_{j=1}^nq_j m^j\frac{d\bar{s}_m}{\bar{s}_ma_m^2} = 0.
\]

Consequently, $\displaystyle q_j = -\sum_{m\in \Sigma(1)\setminus S_i} w_mm^j\frac{a_j^2}{a_m^2}\rho$.
\hfill\rule{2.1mm}{2.1mm}\\

Define $F=a=(a_1,\cdots,a_n): {\cal X} \rightarrow \mathbb{R}^n$. Let $A_i(x) = \displaystyle\min_{m\in \Sigma(1)\setminus S_i} a_m(x)$. For $\eta>0$, consider the domain $U_{i,\eta} = \{x\in U_\eta|A_i(x) \geq A_{i'}(x)\ {\rm for}\ i'\in \Sigma(n)\}$, where $U_\eta = \{x\in {\cal X}|a_m(x) \geq \eta\ {\rm for}\ m\in \Sigma(1)\}$. Notice that proposition \ref{ce} implies that $A_i(x) \geq \lambda_1\tau >0$ for $x \in U_{i,\eta}$. It is easy to observe that there exist $c'>c>0$ such that $[\eta,c\tau]^n \subset F(X_t\cap U_{i,\eta}) \subset [\eta,c'\tau]^n$. For $\omega_t$ and $W$ as in the previous example, we have\\
 
\begin{prop}
\label{ec}
There exists a constant $C_{i,\eta} \geq 0$, such that
\[
\frac{\displaystyle\int_{X_t\cap U_{i,\eta}}\|\bar{\partial} W\|^2 \omega^n_t}{\displaystyle\int_{X_t\cap U_{i,\eta}} \omega^n_t} = \frac{C_{i,\eta}+O(\tau^{-1}\log \tau)}{|\log |t|^2|^3}.
\] 
\end{prop}
{\bf Proof:} 
We may compute the volume of $X_t\cap U_{i,\eta}$.\\
\[
\int_{X_t\cap U_{i,\eta}} \omega^n_t = n!2^n \int_{F(X_t\cap U_{i,\eta})} \rho \prod_{j=1}^n\frac{da_j}{a_j^2} = n!2^n \int_{[\eta,c\tau]^n} \rho \prod_{j=1}^n\frac{da_j}{a_j^2}(1+O(1/\tau))
\]
\[
= n!2^n \prod_{j=1}^n\left(\int_\eta^{c\tau} \frac{da_j}{a_j^2}\right)(1+O(\tau^{-1}\log \tau)) = \frac{n!2^n}{\eta^n}(1+O(\tau^{-1}\log \tau)).
\]

Notice

\[
\bar{\partial}\left(\frac{a_j^2}{a_m^2}\rho\right) = \rho\frac{2a_j}{a_m^2}\bar{\partial}a_j + \frac{a_j^2}{a_m^2}O\left(\frac{\bar{\partial}a_{j'}}{a_m}\right).
\]\\
It is straightforward to compute\\
\[
\left\|\sum_{m\in \Sigma(1)\setminus S_i} w_mm^j\bar{\partial}\left(\frac{a_j^2}{a_m^2}\rho\right)\right\|^2 = 4a_j^4\rho\left|\sum_{m\in \Sigma(1)\setminus S_i}\frac{w_mm^j}{a_m^2}\right|^2.
\]

According to lemma \ref{eb}, we have

\[
\bar{\partial} W = -\sum_{j=1}^n\sum_{m\in \Sigma(1)\setminus S_i} w_mm^j\bar{\partial}\left(\frac{a_j^2}{a_m^2}\rho\right) z_j\frac{\partial}{\partial z_j}
\]
\[
\|\bar{\partial} W\|^2 = \sum_{j=1}^n 4a_j^2\rho\left|\sum_{m\in \Sigma(1)\setminus S_i}\frac{w_mm^j}{a_m^2}\right|^2.
\]
\[
\int_{X_t\cap U_{i,\eta}}\|\bar{\partial} W\|^2 \omega^n_t = n!2^n \int_{F(X_t\cap U_{i,\eta})} \sum_{j=1}^n 4\rho\left|\sum_{m\in \Sigma(1)\setminus S_i}\frac{w_mm^j}{a_m^2}\right|^2 da_j\prod_{j'\not=j}\frac{da_{j'}}{a_{j'}^2}
\]

For each $j$, let

\[
\tilde{U}_{ij,\eta}^0 = \{a\in \mathbb{R}^n|\eta\leq a_j \leq c_j\tau,\ \eta\leq a_{j'} \leq c\tau,\ {\rm for}\ j'\not=j\}.
\]
\[
\tilde{U}_{ij,\eta}^1 = \{a\in F(X_t\cap U_{i,\eta})|\eta\leq a_{j'} \leq c\tau,\ {\rm for}\ j'\not=j\},\ \tilde{U}_{ij,\eta}^2 = F(X_t\cap U_{i,\eta})\setminus \tilde{U}_{ij,\eta}^1.
\]

It is straightforward to derive that 

\[
\int_{\tilde{U}_{ij,\eta}^2} 4\rho\left|\sum_{m\in \Sigma(1)\setminus S_i}\frac{w_mm^j}{a_m^2}\right|^2 da_j\prod_{j'\not=j}\frac{da_{j'}}{a_{j'}^2} = O\left(\frac{1}{\eta^{n-1}\tau^4}\right)
\]
\[
\int_{\tilde{U}_{ij,\eta}^1} 4\rho\left|\sum_{m\in \Sigma(1)\setminus S_i}\frac{w_mm^j}{a_m^2}\right|^2 da_j\prod_{j'\not=j}\frac{da_{j'}}{a_{j'}^2}
\]
\[
=\int_{\tilde{U}_{ij,\eta}^1} 4\rho_j(b_j)\left|\sum_{m\in \Sigma(1)\setminus S_i}\frac{w_mm^j}{(w_m\tau + m^ja_j)^2}\right|^2 da_j\prod_{j'\not=j}\frac{da_{j'}}{a_{j'}^2} + O\left(\frac{\log \tau}{\eta^{n-1}\tau^4}\right)
\]
\[
\left(\int_{\tilde{U}_{ij,\eta}^1}-\int_{\tilde{U}_{ij,\eta}^0}\right) 4\rho_j(b_j)\left|\sum_{m\in \Sigma(1)\setminus S_i}\frac{w_mm^j}{(w_m\tau + m^ja_j)^2}\right|^2 da_j\prod_{j'\not=j}\frac{da_{j'}}{a_{j'}^2} = O\left(\frac{\log \tau}{\eta^{n-1}\tau^4}\right)
\]
\[
\int_{\tilde{U}_{ij,\eta}^0} 4\rho_j(b_j)\left|\sum_{m\in \Sigma(1)\setminus S_i}\frac{w_mm^j}{(w_m\tau + m^ja_j)^2}\right|^2 da_j\prod_{j'\not=j}\frac{da_{j'}}{a_{j'}^2} 
\]
\[
= \frac{4}{\eta^{n-1}|\log |t|^2|^3} \sum_{j=1}^n B_j \prod_{j'\not=j}\int_1^{+\infty} \frac{dx_{j'}}{x_{j'}^2} + O\left(\frac{\log \tau}{\eta^{n-1}\tau^4}\right),
\]

where

\[
B_j = \int_0^{c_j} \rho_j(b_j)\left|\sum_{m\in \Sigma(1)\setminus S_i}\frac{w_mm^j}{(w_m + m^jb_j)^2}\right|^2db_j,
\]

with $c_j = \frac{w_{\tilde{m}_j}}{1-\tilde{m}_j^j}$, $b_j=a_j/\tau$, $x_{j'} = a_{j'}/\eta$, and $\rho_j(b_j)$ is $\rho$ replacing $a_j/a_m$ by $b_j/(w_m + m^jb_j)$ and replacing $a_{j'}$ for $j'\not=j$ by zero. Combining all these estimates, we have

\[
\int_{X_t\cap U_{i,\eta}}\|\bar{\partial} W\|^2 \omega^n_t = \frac{n!2^{n+2}}{\eta^{n-1}}\frac{1}{|\log |t|^2|^3} \sum_{j=1}^n (B_j +O(\tau^{-1}\log \tau)),
\]

We may take $\displaystyle C_{i,\eta} = 4\eta\sum_{j=1}^n B_j$ for the proposition to hold.
\hfill\rule{2.1mm}{2.1mm}\\

\subsection{the toroidal case}

With respect to the local \k metric $\displaystyle\omega_p = \hat{\omega} + \frac{i}{2\pi}\partial\bar{\partial}\log h_p$ and parametrizing function $t$ on $U_p$, we can similarly define $W_{(p)}= \frac{\nabla \log t}{|\nabla \log t|^2}$. Let $\displaystyle W = \sum_{p\in \Sigma}\mu_pW_{(p)}$. $\bar{\partial} W$ also represents the Kodaira-Spencer deformation class. We have

\begin{prop}
\label{ea}
There exists a constant $C>0$ independent of $t$ such that

\[
\int_{X_t}\|\bar{\partial} W\|_{g_t}^2 \omega^n_t \leq \frac{C}{|\log |t|^2|^3}\int_{X_t} \omega^n_t.
\]
\end{prop}
{\bf Proof:} 
Locally in each $U_{q\sigma}^0$, we will use similar coordinate and proper basis $\{W_j,\bar{W}_j\}_{j=1}^n$ as in the proof of proposition \ref{cc}. Then the dual basis is $\{\beta_j,\bar{\beta}_j\}_{j=1}^n$, where $\beta_j = \frac{dz_j}{a_jz_j}$, $a_j = \log |z_j|^2$ for $1\leq j \leq l$ and $\beta_i = dz_i$ for $l+1 \leq j \leq n$. Recall that $O(1)$ denotes a smooth function on terms in (\ref{ca}) and other smooth bounded terms. (Notice that here we assume $a_m = \log|s_m|^2$, which is slightly different from (\ref{ca}) and do not affect our arguments here. In this proof, we are using $a_j$ to denote $a_{m_j}$ and $O(a_j/a_m)$ as a shorthand for $O(a_j/a_m,1\leq j \leq n,m\in \Sigma_q(1)\setminus S_\sigma)$.) We will also use $O(1)$ to denote a tensor with $O(1)$ coefficients with respect to the proper and dual proper basis. It is easy to see that the action of the proper basis $\{W_j,\bar{W}_j\}_{j=1}^n$ will send $O(1)$ to $O(1)$, also $\bar{\partial} W_j = O(1)$. Under such notation, we have

\[
\omega_{t,q} = \sum_{j,k=1}^n g_{j\bar{k}}\beta_j\bar{\beta}_k.
\]
\[
g_{j\bar{k}} = \delta_{jk}\left(1+\frac{1}{a_k}O(1)\right) + a_ja_kO\left(\frac{1}{a_m^2}\right) + \frac{1}{a_ja_k}O(1),\ {\rm for}\ 1\leq j,k \leq l.
\]
\[
g_{j\bar{k}} = \frac{1}{a_k}O(1),\ g^{j\bar{k}} = \frac{1}{a_k}O(1),\ {\rm for}\ 1\leq k \leq l\ {\rm and}\ l+1 \leq j\leq n.
\]

It is straightforward to derive that

\[
\left.i\left(t\frac{\partial}{\partial t}\right)\omega_q\right|_{X_t} = \sum_{j=1}^l\sum_{m\in \Sigma_q(1)\setminus S_\sigma} w_mm^j\frac{a_j}{a_m^2}\bar{\beta}_j + O\left(\frac{1}{a_m^2}\right),
\]
\[
W_{(q)} = t\frac{\partial}{\partial t} - \sum_{j=1}^l\sum_{m\in \Sigma_q(1)\setminus S_\sigma} w_mm^j\frac{a_j}{a_m^2}\rho W_j + O\left(\frac{1}{a_m^2}\right),
\]
\[
\bar{\partial} W_{(q)} = -\sum_{j=1}^l\sum_{m\in \Sigma_q(1)\setminus S_\sigma} w_mm^j\bar{\partial}\left(\frac{a_j^2}{a_m^2}\rho\right) z_j\frac{\partial}{\partial z_j}
 + O\left(\frac{1}{a_m^2}\right).
\]

Applying proposition \ref{ec}, we can find $C>0$ independent of $t$ such that 

\[
\int_{U_{q\sigma}^0\cap X_t}\|\bar{\partial} W_{(q)}\|_{g_t}^2 \omega^n_t \leq \frac{C}{|\log |t|^2|^3}\int_{U_{q\sigma}^0\cap X_t} \omega^n_t.
\]

The rest of the proof closely resemble the proof of lemma \ref{ch}. For any $x\in X_t$, there exist a $q\in \Sigma$ such that $x\in X_t\cap U_q^0$. Since

\[
\sum_{p\in \Sigma} \mu_p =1,\ \ \sum_{p\in \Sigma} \bar{\partial}\mu_p =0.
\]

We have

\[
\sum_{p\in \Sigma}\bar{\partial} \mu_pW_{(p)} = \sum_{p\in \Sigma}\bar{\partial} \mu_p(W_{(p)} - W_{(q)}).
\]

Since $U_q^0 \cap {\rm supp}(\mu_p)=\emptyset$ when $D_q\not\subset \bar{D}_p$ according to (\ref{ad}), we may consider only those $p\in \Sigma$ satisfying $D_q\subset \bar{D}_p$. As in the proof of lemma \ref{ch}, for such $p,q\in \Sigma$, we can naturally define $\Sigma_q^q \cap\Sigma_p^p$. For any $m\in\Sigma_p^p \setminus \Sigma_q^q$, $(s_m)\cap \bar{D}_q=\emptyset$. Consequently, $\|s_m\|_m^2= 1$ and $a_m = \eta$ on $\tilde{U}_q$ for $m\in\Sigma_p^p \setminus \Sigma_q^q$. Hence

\[
W_{(p)} - W_{(q)} = \sum_{j=1}^lO\left(\frac{a_j}{a_m^2}\right)W_j + O\left(\frac{1}{a_m^2}\right)
\]
\[
\bar{\partial}W_{(p)} = \sum_{j=1}^lO\left(\frac{a_j}{a_m^2}\right) + O\left(\frac{1}{a_m^2}\right)
\]

on ${\rm supp}(\mu_p)\cap U_q \subset U_p\cap U_q$. From the explicit expressions of $\mu_p$, it is straightforward to check that $\bar{\partial}\mu_p = O(1/\log \tau)$ with respect to the Hermitian metric $\omega_t$. Consequently 

\[
\int_{U_{q\sigma}^0\cap X_t}\left\|\sum_{p\in \Sigma} \bar{\partial}\mu_pW_{(p)}\right\|_{g_t}^2 \omega^n_t \leq \frac{C}{\tau^3\log\tau}\int_{U_{q\sigma}^0\cap X_t} \omega^n_t.
\]

\[
\int_{U_{q\sigma}^0\cap X_t}\left\|\sum_{p\in \Sigma} \mu_p\bar{\partial} W_{(p)}\right\|_{g_t}^2 \omega^n_t \leq \frac{C}{\tau^3}\int_{U_{q\sigma}^0\cap X_t} \omega^n_t.
\]

Combine these estimates for all $\sigma\in \Sigma_q(l)$, $q\in \Sigma$ applying to

\[
\bar{\partial} W = \sum_{p\in \Sigma}\mu_p\bar{\partial} W_{(p)} + \sum_{p\in \Sigma}\bar{\partial} \mu_pW_{(p)},
\]

we will get the desired estimate.
\hfill\rule{2.1mm}{2.1mm}\\

{\bf Remark:} It is not hard to observe that the constant $C_{i,\eta} \geq 0$ in proposition \ref{ec} is actually positive. With this observation and a bit more argument, one can show that the lower bound estimate in proposition \ref{ea} (more precisely the estimate in proposition \ref{ea} with the reversed inequality) is also true. Since such more precise estimates are not needed for arguments in this paper, we will omit them here.\\

{\bf Proof of theorem \ref{ab}:} As pointed out in \cite{Tian1},

\[
g_{WP}\left.\left(\frac{d}{dt},\frac{d}{dt}\right)\right|_{X_t} = \int_{X_t}\left\|H\left(\frac{d}{dt}\right)\right\|^2_{g_{E,t}} \omega_{E,t}^n,
\]

where $\displaystyle H\left(\frac{d}{dt}\right)$ denote the harmonic representative of the Kodaira-Spencer deformation class. As mentioned earlier, such class can also be represented by $\displaystyle \frac{\bar{\partial} W}{t}$. Applying proposition \ref{ea} and theorem \ref{da}, we have

\[
\int_{X_t}\left\|H\left(\frac{d}{dt}\right)\right\|^2_{g_{E,t}} \omega_{E,t}^n \leq \int_{X_t}\left\|\frac{\bar{\partial} W}{t}\right\|^2_{g_{E,t}} \omega_{E,t}^n \leq C\int_{X_t}\left\|\frac{\bar{\partial} W}{t}\right\|^2_{g_t} \omega_t^n \leq \frac{C}{|\log |t||^3|t|^2}.
\]
\hfill\rule{2.1mm}{2.1mm}\\

\ifx\undefined\bysame
\newcommand{\bysame}{\leavevmode\hbox to3em{\hrulefill}\,}
\fi


\noindent
\begin{thebibliography}{100}

\bibitem{A} T. Aubin, {\em Equation du type de \ma sur les vari\'{e}t\'{e}s K\"{a}hleriennes compacts}, C. R. Acad. Sci. Prais {\bf 283} (1976), 119-121.


\bibitem{CG} J. Cheeger and M. Gromov, {\em Collapsing Riemannian manifolds while keeping their curvature bounded. I}, J. Differential Geom. {\bf 23} (1986), 309-346; II, J. Differential Geom. {\bf 32} (1990), 269-298.


\bibitem{CY}
S. Y. Cheng and S.-T. Yau, {\em On Inequality between Chern numbers of singular \k surfaces and characterization of orbit space of discrete group of $SU(2,1)$}, Contemporary Math. {\bf 49} (1986), 31-43.

\bibitem{Lu}
N. Leung and P. Lu, {\em Degeneration of \k Einstein metrics on complete \k manifolds}, Comm. Analysis and Geom. {\bf 7} (1999), 431-449. 

\bibitem{Mu}
G. Kempf, F. Knudsen, D. Mumford and B. Saint-Donat, {\em Toroidal Embeddings I}, Lecture Notes in Mathematics {\bf 339}, Springer-Verlag 1973. 

\bibitem{ruan1} 
W.~D. Ruan, {\em On the convergence and collapsing of K\"{a}hler manifolds}, Journal of Differential Geometry, Volume 52 (1999), 1-40.

\bibitem{ruan2}
\bysame, {\em Canonical coordinates and Bergmann metrics}, Communications in Analysis and Geometry, Vol. 6 (1998), 589-631.
%
\bibitem{ruan3}
\bysame, {\em Degeneration of \ke manifolds I: The normal crossing case}, To appear in Comm. Contemporary Math.

\bibitem{Tian1} 
G. Tian, {\em Degeneration of \ke manifolds I}, Proceedings of Symposia in Pure Mathematics, Vol.54, Part 2, 595-609.

\bibitem{TianYau} 
G. Tian and S.-T. Yau, {\em existence of \ke metrics on complete \k manifolds and their applications to algebraic geometry}, Math. Aspects of String Theory (Edited by S. -T. Yau) pp. 574-628, World Sci. Publishing, 1987.

\bibitem{yau}
S.-T. Yau, {\em On {C}alabi's conjecture and some new results in algebraic
  geometry}, Proc. Nat. Acad. Sci. U.S.A. {\bf 74} (1977), 1798--1799.

\bibitem{yau1}
\bysame, {\em M\'{e}triques de \ke sur les vari\'{e}t\'{e}s overtes}, Ast\'{e}risque {\bf 58} (1978), 163-167.

\bibitem{Yau2}
\bysame, {\em On the Ricci curvature of a compact K\"{a}hler manifold and the complex Monge-Amp\`{e}re equation, I}, Comm. Pure. and Appl. Math.,  {\bf 31} (1978), 339-411.
\end{thebibliography}
\end{document}